\documentclass[a4paper,reqno,12pt]{amsart}

\usepackage{amssymb}
\usepackage{amsthm}

\newtheorem{theorem}{Theorem}[section]
\newtheorem{lemma}[theorem]{Lemma}

\numberwithin{equation}{section}

\begin{document}

\author[Aimo Hinkkanen]{Aimo Hinkkanen}
\address{Department of Mathematics, University of Illinois Urbana-Champaign,
 Urbana, IL 61801 USA }
\email{aimo@illinois.edu}
\author[Ilpo Laine]{Ilpo Laine}
\address{Department of Physics and Mathematics, University of Eastern Finland,
Joensuu, 80100 Finland}
\email{ilpo.laine@uef.fi}

\dedicatory{Dedicated to the memory of Lawrence Zalcman}



\title[Value sharing and Stirling numbers]{Value sharing and Stirling numbers}

\begin{abstract}
Let $f$ be an entire function and $L(f)$ a linear differential polynomial in $f$ with constant coefficients. Suppose that $f$, $f'$, and $L(f)$ share a meromorphic function $\alpha(z)$ that is a small function with respect to $f$. A characterization of the possibilities that may arise was recently obtained by Lahiri. However, one case leaves open many possibilities.
We show that this case has more structure than might have been expected, and that a more detailed study of this case involves, among other things, Stirling numbers of the first and second kinds. We prove that the function $\alpha$ must satisfy a linear homogeneous differential equation with specific coefficients involving only three free parameters, and then $f$ can be obtained from each solution. Examples suggest that only rarely do single-valued solutions $\alpha(z)$ exist, and even then they are not always small functions for $f$.
\end{abstract}

\maketitle

\section{Introduction}\label{intro}

We show that one case of a value sharing problem investigated by Lahiri in \cite{Lah} has more structure than might have been expected, and that a more detailed study of this case involves, among other things, Stirling numbers of the first and second kinds.

We assume that the reader is familiar with Nevanlinna's value distribution theory (see, for example, \cite{Hay}). Unless otherwise mentioned, all functions are assumed to be meromorphic functions in the complex plane ${\mathbb C}$. We say that functions $F$ and $G$ share a function (possibly constant) $a$ counting multiplicities (CM) if $F-a$ and $G-a$ have the exactly same zeros in ${\mathbb C}$, and the multiplicities of the zeros are the same at each zero.

Our starting point is formed by the papers \cite{WL} due to Wang and Laine and \cite{Z} by Zhang.
For the convenience to the reader, we recall the following theorem from \cite{Z}:

\begin{theorem}\label{thm A} Suppose that $f$ is a non-constant entire function, $\alpha (z)$ is a non-zero small function with respect to $f$, and $n\geq 2$ is an integer. Let $L(f):=a_{n}f^{(n)}+a_{n-1}f^{(n-1}+\cdots +a_{1}f'+a_{0}f$ be a linear differential polynomial to $f$ with constant coefficients such that $a_{0}\neq 0$. If $f,f'$ and $L(f)$ share $\alpha (z)$ CM, then one of the following assertions holds:

(i) $f(z)=be^{z}$, where $b$ is a non-zero constant.

(ii) $\alpha (z)$ reduces to a constant, say $a$, $a_{0}+\sum_{j=1}^{n}a_{j}c^{j-1}=1$ for some constant $c\neq 0$, and
$$f(z)=be^{cz}-\frac{a(1-c)}{c},$$
where $b\neq 0$ is a constant.

(iii) $f(z)=\lambda e^{cz}+p(z)+a(z),$ where $c$ is a constant, $p(z)$ and $a(z)$ are two polynomials satisfying $p'-cp=a-a'$, and $\deg a\leq n-3.$
\end{theorem}

The case (iii) here is a completion to Theorem 1.1 in \cite{WL}. However, as pointed out in a recent paper \cite{Lah} by Lahiri, even the above theorem due to Zhang is not complete. The key result of Lahiri is \cite{Lah}, Theorem 1.1, completing these two preceding papers. For the convenience of the reader, we recall this result due to Lahiri:

\begin{theorem}\label{Thm B} Let $f$ be a non-constant entire function, $\alpha$ be a small function relative to $f$, $n\geq 2$ be a positive integer and $L(f):=a_{0}f+a_{1}f'+\cdots +a_{n}f^{(n)}$ be with constant coefficients, $a_{n}\neq 0$. Suppose that $f$ and $f'$ share $\alpha$ CM and $f$ and $L(f)$ share $\alpha$ CM as well, except possibly in a finite set. Then five possibilities may appear:

\medskip

(i) $f(z)= \alpha +\exp ((\lambda e^{cz}+\mu )/c)$ and $L(f)= \alpha +a_{n}\lambda^{n}\exp ((\lambda e^{cz}+nc^{2}z+\mu )/c)$, where $ \alpha'\neq  \alpha $ and $c\neq 0, \lambda\neq 0$ and $\mu$ are constants.

\medskip

(ii) $f(z)=\lambda e^{z}$ and $L(f)=f$, where $ \alpha'\neq  \alpha $ and $\lambda\neq 0$ is a constant.

\medskip

(iii) $f(z)=\lambda e^{cz}+ \alpha +p$ and $L(f)=((1-c)a_{0}+c)f+(1-c)(1-a_{0}) \alpha $, where $ \alpha $ is a non-constant polynomial and $\lambda\neq 0$, $c\neq 0,1$ are constants. Moreover, $p=\frac{c-1}{c}\left( -\frac{ \alpha }{c-1}+\sum_{j=1}^{u}\frac{ \alpha ^{(j)}}{c^{j}}\right)$ is a polynomial, $u=\deg  \alpha $ and $p'-cp= \alpha - \alpha'$.

\medskip

(iv) $f(z)=\lambda e^{cz}+(1-1/c) \alpha $ and $L(f)= \alpha +((1-c)a_{0}+c)(f- \alpha )$, where $ \alpha $ is a non-zero constant and $\lambda\neq 0, c\neq 0,1$ are constants.

\medskip

(v) $f$ satisfies the differential equations $f'=\lambda e^{cz}f+(1-\lambda e^{cz})\alpha$ and $L(f)=a_{n}\lambda^{n}e^{ncz}f+(1-a_{n}\lambda^{n}e^{ncz})\alpha$, where $\alpha$ is non-constant, $\alpha '\neq\alpha$, and $c,\lambda$ are non-zero constants. Moreover, $a_{0}=0$ and $a_{n-1}\neq 0$. If $n=2$, then $a_{1}+a_{2}c=0$, and if $n\geq 3$, then $\sum_{k=1}^{n}(-1)^{k-1}a_{k}\left(\frac{2a_{n-1}}{n(n-1)a_{n}}\right)^{k-1}=0$.
\end{theorem}

In this paper, we   analyze the case (v) of the above theorem in more detail, trying to demonstrate that the fifth case actually may appear in some special cases only.

We denote the Stirling numbers of the first kind by $s(n,k)$ and those of the second kind by $S(n,k)$. We will recall the definitions and basic properties of these numbers in Section~\ref{Stir}.  We obtain the following result.

\begin{theorem} \label{th0}
Suppose that $f$, $n$, $L(f)$, and $\alpha$ are as in case (v) of Theorem~\ref{Thm B}, so that $(1-\lambda e^{cz})\alpha$  is an entire function. Then $a_0=0$ and
\begin{equation} \label{eq35ee}
a_j =
a_n c^{n-j} s(n,j)
\end{equation}
for $1\leq j\leq n$.
Furthermore, $\alpha(z)$ satisfies the linear homogeneous differential equation
\begin{eqnarray} \label{a2a}
0 & = &
\left(1-
a_n \sum_{p=0}^{n-1} c^{n - p - 1} \lambda^p e^{pcz}
(-1)^{n - p - 1} (n - p - 1)!  \right) \alpha \notag
\\
& {} &
-  a_n \sum_{k=1}^{n-2}   \alpha^{(k)}
c^{n-1-k} \left( s(n,k+1) + \sum_{p=1}^{n-k} (\lambda/c)^p e^{pcz}
(s(n - p,k+1) - c s(n - p,k) ) \right) \notag
\\
& {} &
-a_n ( 1 -  \lambda  e^{cz}           ) \alpha^{(n-1)}
\end{eqnarray}
 of order $n-1$. When $n=2$, the middle line in (\ref{a2a}) represents an empty sum, which we understand to be equal to zero. Finally, $f$ can be solved for from the equation
\begin{equation}\label{eq1a}
f'(z)=\lambda e^{cz}f(z)+(1-\lambda e^{cz})\alpha (z)  .
\end{equation}
\end{theorem}

Theorem~\ref{th0} leaves open various questions. If $n\geq 2$ is given, for which choices of the non-zero complex parameters $\lambda$, $c$, and $a_n$, does the equation (\ref{a2a}) have a non-trivial solution $\alpha(z)$ that is meromorphic in the complex plane? When such a solution $\alpha$ exists, when is the function $f$ obtained from (\ref{eq1a}) entire, and if it is, when is $\alpha$ a small function for $f$? We will discuss these questions for $n=2$ and $n=3$. The general case remains open.

\medskip

This paper has been organized as follows. After some elementary remarks (mainly for the case $n=2$) in Section \ref{initial}, we recall some basic facts of Stirling numbers in Section \ref{Stir}. These facts will then be extensively used in the next Section \ref{derf} to construct expressions for the derivatives $f^{(n)}$ of $f$. The next two sections will then be devoted to the differential equation satisfied $\alpha (z)$, first constructing the equation itself in Section \ref{alfa}, and then considering the possibilities to solve the equation in Section \ref{alfa2}.

\section{Initial remarks on Lahiri's work}\label{initial}

We first proceed to considering the case $n=2$. In this case, the differential equations satisfied by $f$ may be written as follows:
\begin{equation}\label{eq1}
f'(z)=\lambda e^{cz}f(z)+(1-\lambda e^{cz})\alpha (z)
\end{equation}
and
\begin{equation}\label{eq2}
L(f)(z)=a_{2}\lambda^{2}e^{2cz}f(z)+(1-a_{2}\lambda^{2}e^{2cz})\alpha (z).
\end{equation}

The solution to equation (\ref{eq1}) may  be written as
\begin{equation}\label{sol1}
f(z)=C\exp ({(\lambda /c)e^{cz}})+\exp ({(\lambda /c)e^{cz}})\int\exp ({-(\lambda /c)e^{cz}})(1-\lambda e^{cz})\alpha (z)dz,
\end{equation}
$C$ being a constant of integration. Obviously, we now have
$$\frac{f'-\alpha}{f-\alpha}=\lambda e^{cz},$$
so the CM-condition for $f,f'$ is satisfied. Observe that $(1-\lambda e^{cz})\alpha (z)$ in the expression of $f(z)$ is entire by (\ref{eq1}), hence the integration in (\ref{sol1}) brings no problems (although it usually cannot be expressed in explicit form). To verify the CM-condition for $f$ and $L(f)$, we may proceed by recalling $a_{1}=-a_{2}c$ from Theorem \ref{Thm B}, differentiating (\ref{eq1}), and substituting this and (\ref{eq1}) into
$$-a_{2}cf'(z)+a_{2}f''(z)-a_{2}\lambda^{2}e^{2cz}f(z)=(1-a_{2}\lambda^{2}e^{2cz})\alpha (z).$$

We  obtain, after simplification, that
\begin{equation}\label{eq3}
a_{2}(1-\lambda e^{cz})\alpha '(z)- (1+a_{2}c-a_{2}\lambda e^{cz})\alpha (z)=0.
\end{equation}
It is now immediate to verify that
\begin{equation}\label{sol4}
\alpha (z)=\widetilde{C}e^{(1+\frac{1}{a_{2}c})cz}(\lambda e^{cz}-1)^{-1+\frac{1}{c}-\frac{1}{a_{2}c}}.
\end{equation}

Recalling that
$$(1-\lambda e^{cz})\alpha (z)=-\widetilde{C}e^{(1+\frac{1}{a_{2}c})cz}(\lambda e^{cz}-1)^{\frac{1}{c}-\frac{1}{a_{2}c}}$$
has to be an entire function, we must have
\begin{equation}\label{enta2}
\frac{1}{c}-\frac{1}{a_{2}c}=s\in\{ 0\}\cup\mathbb{N}.
\end{equation}

Moreover, we must have $sc\neq 1$. Therefore, $a_{2}=1/(1-sc)$, where $s\in\{ 0\}\cup\mathbb{N}$. In particular, if $s=0$, then $a_{2}=1$.

\medskip

\textbf{Example.} Take  $s=1$. Then $a_{2}=\frac{1}{1-c}$ and $a_{1}=-a_{2}c=\frac{c}{c-1}$. In this special case, see (\ref{sol4}), we get $\alpha (z)=e^{z}$, with the constant of integration being $= 1$. Recalling (\ref{sol1}), we obtain
\begin{equation}\label{fs1}
f(z)=\exp \left(  \frac{\lambda}{c}e^{cz}  \right)+e^{z}.
\end{equation}
Computing $L(f)=\frac{c}{c-1}f'(z)+\frac{1}{1-c}f''(z)$, we get
$$L(f)=e^{z}+\frac{1}{1-c}\lambda^{2}e^{2cz}\exp \left(\frac{\lambda}{c}e^{cz} \right).$$
By computation, we now get
$$(f'-\alpha )/(f-\alpha )=(f'-e^{z})/(f-e^{z})=\lambda e^{cz}$$
and
$$(L(f)-\alpha )/(f-\alpha )=(L(f)-e^{z})/(f-e^{z})=\frac{1}{1-c}\lambda^{2}e^{2cz}=a_{2}\lambda^{2}e^{2cz}.$$
Observe that in this case, we clearly have that $\alpha (z)=e^{z}$ is a small function with respect to $f$ in (\ref{fs1}).

\medskip

We next proceed to considering the cases when $s=0$ or $s\geq 2$. Then the integral in (\ref{sol1}) cannot be expressed in explicit form. Indeed, suppose that (\ref{sol4}) holds where $\frac{1}{c}-\frac{1}{a_{2}c}=s$ is a non-negative integer. Then with the notation $u= \lambda e^{cz}$ the integrand on the right hand side of (\ref{sol1}) is a sum where $m$ runs from $0$ to $s$ and where the generic term is a constant multiple (with the constant depending on $m$) of
$$
 e^{-(1/c)u} u^{m-1-s+(1/c)} \, du  .
$$
With $v=-(1/c)u$, this is a constant multiple of
$$
e^{v} v^{m-1-s+(1/c)} \, dv  .
$$
This can be integrated explicitly if $m-1-s+(1/c)$ is a non-negative integer, and hence for every $m$ with $0\leq m\leq s$ if $1/c$ is an integer with $1/c\geq s+1$.
Assuming that $1/c=\nu\geq s+1$ we then have
$a_2 = \nu/(\nu-s)$ so that $a_2$ is a rational number with $a_2\geq 1$. In other cases the integral can be expressed in terms of the incomplete gamma function.

\medskip

\textbf{Remark.} We next point out a couple of necessary conditions for the existence of an entire function $f$ to satisfy equations (\ref{eq1}) and (\ref{eq2}). Suppose first that $\alpha$ is entire. Looking at points $\tilde{z}$ such that $\lambda e^{c\tilde{z}}=1$, then it is immediate to observe that $f'(\tilde{z})-f(\tilde{z})=0$. To consider the case when $\alpha$ is meromorphic in general, eliminating $\alpha$ from equations (\ref{eq1}) and (\ref{eq2}), we obtain

\begin{eqnarray}\label{eq4}
&{}&
\frac{\lambda e^{cz}(1-a_{n}\lambda^{n-1}e^{(n-1)cz})f+(a_{1}(1-\lambda e^{cz})-(1-a_{n}\lambda^{n}e^{ncz}))f'}{1-\lambda e^{cz}}
\notag
\\
&{}&
=-(a_{2}f''+\cdots +a_{n}f^{(n)}).
\end{eqnarray}

The right hand side of (\ref{eq4}) is entire, hence the same applies with the left hand side as well. This means that the numerator
$$\lambda e^{cz}(1-a_{n}\lambda^{n-1}e^{(n-1)cz})f+(a_{1}(1-\lambda e^{cz})-(1-a_{n}\lambda^{n}e^{ncz}))f'$$
vanishes at all points $\tilde{z}$ such that $\lambda e^{c\tilde{z}}=1$. This means that we get two possible necessary conditions: Either (1) $a_{n}=1$, or (2) $f'(\tilde{z})-f(\tilde{z})=0$.

\section{Stirling numbers}\label{Stir}

In what follows, it appears that Stirling numbers may be used to obtain a formula for the derivatives $f^{(n)}$ in terms of $f$. Indeed, from (\ref{eq1}) we can get a formula for $f^{(n)}$ of the form
\begin{equation} \label{eq5}
f^{(n)} = A_n f + B_n
\end{equation}
where neither $A_n$ nor $B_n$ involves $f$ or its derivatives. We will find explicit expressions for $A_n$ and $B_n$. They will involve Stirling numbers of the second kind, see, e.g., \cite{Cha}, Chapter~8.

Stirling numbers of the second kind $S(n,k)$ are defined as the number of ways of partitioning a set of $n$ elements into $k$ (non-empty) subsets. For this to make sense, we must have $1\leq k\leq n$. In order to use $S(n,k)$ in formulas in a more general way, it is customary to define $S(n,k)=0$ if $1\leq n<k$; to define $S(n,0)=0$ if $n\geq 1$; and to define $S(0,0)=1$. We also define $S(n,-1)=0$ for $n\geq 1$.

We will need the identity
\begin{equation} \label{eq6}
S(n,k) = S(n-1,k-1) + k S(n-1,k) .
\end{equation}
To prove this, one can argue as follows. To partition $\{1,2,\dots ,n\}$ into $k$ subsets, we can partition $\{1,2,\dots ,n-1\}$ into $k-1$ subsets and use $\{n\}$ as the $k^{ {\rm th} }$ subset. This gives rise to the term $S(n-1,k-1)$. Or we can partition $\{1,2,\dots ,n-1\}$ into $k$ subsets and for each choice, include the element $n$ in any one of those $k$ subsets. This gives rise to the term $k S(n-1,k)$.

There are many formulas for $S(n,k)$. We will not need all of them but we mention some of them for completeness. We have
\begin{equation} \label{eq7}
S(n,k) = \sum  1^{a_1-1} 2^{a_2-1} 3^{a_3-1} \cdots k^{a_k-1}
\end{equation}
where the sum is taken over all positive integers $a_i$ such that $a_1+a_2+\cdots +a_k=n$. Since each $a_i\geq 1$, this obviously requires that $k\leq n$.
We also have
\begin{equation} \label{eq8}
S(n,k) = \frac{1} {k!} \sum_{m=1}^k  (-1)^{k-m} \binom{k}{m}  m^n=  \sum_{m=1}^k \frac{ (-1)^{k-m} } { m! \, (k-m)!    } m^n  .
\end{equation}
We further have
\begin{equation} \label{eq9}
S(n,k) = \frac{ n! } { k! } \sum  \frac{ 1 } {  a_1! a_2! \cdots a_k!  }
\end{equation}
where the sum is taken over all non-negative integers $a_i$ such that $a_1+a_2+\cdots +a_k=n$. And we have
\begin{equation} \label{eq10}
S(n,k) =  \sum \frac{ n! } { a_1! a_2! \cdots a_n! } \left( \frac{1}{ 1! } \right)^{a_1} \left( \frac{1}{ 2! } \right)^{a_2}  \cdots \left( \frac{1}{ n! } \right)^{a_n}
\end{equation}
where the sum is taken over all non-negative integers $a_i$ such that $a_1+a_2+\cdots +a_n=k$ and $a_1+2a_2+3a_3+\cdots +na_n=n$.

\medskip

To explain the connection to Stirling numbers $s(n,k)$ of the first kind,  we begin with the definition. With the notation $(t)_0=1$ and, for $n\geq 1$,
$$
(t)_n = t(t-1)(t-2)\cdots (t-n+1)
$$
we have, for $n\geq 0$,
\begin{equation} \label{eq35cc}
(t)_n = \sum_{k=0}^n s(n,k) t^k
\end{equation}
and this defines $s(n,k)$. Clearly $s(n,0)=0$, $s(n,1)=(-1)^{n-1} (n-1)!$,  and $s(n,n)=1$ for all $n\geq 1$. We have already defined the numbers $S(n,k)$ but they could also be defined for $n\geq 0$ by
$$
t^n= \sum_{k=0}^n S(n,k) (t)_k .
$$
We define $s(n,k)=S(n,k)=0$ when $k>n$ and $s(0,0)=S(0,0)=1$.

We recall some known results from \cite{Cha}, pp.~278, 281. We have
$$
|s(n,k)| = (-1)^{n-k} s(n,k) .
$$
Using the Kronecker symbol $\delta_{m,n}$, defined by $\delta_{n,n}=1$ and by $\delta_{m,n}=0$ for $m\not= n$, we have
\begin{equation} \label{eq35d}
\sum_{r=k}^n s(n,r) S(r,k) = \delta_{n,k} ,
\end{equation}
$$
\sum_{r=k}^n S(n,r) s(r,k) = \delta_{n,k}  .
$$
There is the recursion formula (\cite{Cha}, p.~293)
\begin{equation} \label{z0}
s(n+1,k) = s(n,k-1) - n s(n,k)  .
\end{equation}
By \cite{Cha}, p.~323, we have
\begin{equation} \label{z00}
s(n,n-1) = -\frac{n(n-1)}{2} .
\end{equation}
Moreover, by \cite{Cha}, Corollary 8.2, $s(n,n)=1$.

\medskip

In what follows, we often make use of the well-known identity (from Pascal's triangle)
\begin{equation} \label{bin}
\binom{n}{k} =\binom{n-1}{k} +\binom{n-1}{k-1} .
\end{equation}
We recall that $\binom{n}{k} =0$ if $k<0$ or if $n<k$.

\section{The formula for the derivatives of $f$}\label{derf}

\begin{lemma} \label{le1}
Suppose that (\ref{eq1}) holds. Then for each $n\geq 1$,  (\ref{eq5}) is valid with $A_n$ and $B_n$ given as follows. We have
\begin{equation} \label{eq11}
A_n =  \sum_{k=1}^n b_{n,k} e^{kcz}
\end{equation}
where
\begin{equation} \label{eq12}
b_{n,k}  =  S(n,k) \lambda^k c^{n-k}
\end{equation}
and
\begin{equation} \label{eq13}
B_n =  \sum_{k=0}^{n-1} \beta_{n,k} \alpha^{(k)}
\end{equation}
where
\begin{equation} \label{eq14}
\beta_{n,k}  = \sum_{j=0}^{n-k}  (\zeta_{n,k,j} - c\,  \varepsilon_{n,k,j} )  c^{n - k - j - 1} \lambda^j e^{jcz} .
\end{equation}
Here for $n\geq 1$, $0\leq k\leq n-1$, $0\leq j\leq n-k$, we have
\begin{equation} \label{eq15}
\zeta_{n,k,j}   =   \sum_{m=0}^{n - 1 - k - j} j^m \binom{k+m}{k} S(n-1-k-m,j)
\end{equation}
for $2\leq j\leq n-k-1$, $\zeta_{n,k,n-k} =0$ (hence $\zeta_{n,0,n} =0$), and
\begin{equation} \label{eq16}
\varepsilon_{n,k,j}   =    \sum_{m=0}^{n  - k - j} j^m \binom{k+m}{k} S(n-1-k-m,j-1)
\end{equation}
for $2\leq j\leq n-k-1$, $\varepsilon_{n,k,n-k} =1$  (hence $\varepsilon_{n,0,n} =1$).
The equations above are valid when $j\geq 2$, while
\begin{equation} \label{eq16a}
\zeta_{n,k,1}   =   \sum_{m=0}^{n  - k - 2}  \binom{k+m}{k} S(n-1-k-m,1)
=  \sum_{m=0}^{n  - k - 2}  \binom{k+m}{k} =   \binom{n-1}{k+1}
\end{equation}
when $0\leq k\leq n-2$,  $\zeta_{n,n-1,1} = \binom{n-1}{n}=0$,
and
\begin{equation} \label{eq16c}
\zeta_{n,k,0}   =  0
\end{equation}
when $0\leq k\leq n-2$, $\zeta_{n,n-1,0}=1$ and in particular $\zeta_{1,0,0}=1$,
\begin{equation} \label{eq16d}
\varepsilon_{n,k,1}   =    \binom{n-1}{k}  ,
\end{equation}
and
\begin{equation} \label{eq16e}
\varepsilon_{n,k,0}   =   0  .
\end{equation}
\end{lemma}

Note that if we were to apply  (\ref{eq15}) with $j=n-k$, we would get a sum where $m$ goes from $0$ to $-1$. Interpreting this as an empty sum, hence equal to zero, we also obtain  $\zeta_{n,k,n-k} =0$, which we noted separately. For completeness, we also define $\zeta_{0,k,j}=0$ and $\varepsilon_{0,k,j}=0$ for all $k$ and $j$.

Before proceeding, we obtain some identities for the quantities $\zeta_{n,k,j} $ and $\varepsilon_{n,k,j} $ defined in Lemma~\ref{le1}.

\begin{lemma} \label{le2}
If $n\geq 1$, $1\leq k\leq n$, $0\leq j\leq n+1-k$, we have
\begin{equation} \label{eq16f}
\zeta_{n+1,k,j}   =  j \zeta_{n,k,j}  + \zeta_{n,k-1,j}
\end{equation}
and
\begin{equation} \label{eq16g}
\varepsilon_{n+1,k,j}   =  j \varepsilon_{n,k,j}  + \varepsilon_{n,k-1,j}
\end{equation}
with the convention that $\zeta_{n,k,j}=0$ and $\varepsilon_{n,k,j}=0$ if $k\geq n$ or if $j\geq n-k+1$.
We further have
\begin{equation} \label{eq16h}
\zeta_{1,0,0} = 1, \,\, \varepsilon_{1,0,0} =0, \,\, \zeta_{1,0,1}=0,\,\, \varepsilon_{1,0,1}=1
\end{equation}
and
\begin{equation} \label{eq16i}
\zeta_{n+1,0,j}   =  j \zeta_{n,0,j}  + S(n,j) , \qquad  \varepsilon_{n+1,0,j}   =  j \varepsilon_{n,0,j}  + S(n,j-1)  .
\end{equation}
\end{lemma}

{\bf Proof of Lemma~\ref{le2}.}
The equation (\ref{eq16h}) follows directly from the definitions
(\ref{eq16a})--(\ref{eq16e}).

For $j=0$, (\ref{eq16g}) follows from (\ref{eq16e}) since all quantities vanish.

For $j=0$, (\ref{eq16f}) follows from (\ref{eq16a}) since it implies that
$$
\zeta_{n+1,k,0}   - \zeta_{n,k-1,0}  = 0-0=0
$$
if $1\leq k\leq n-1$, while if $k=n$ then the line after (\ref{eq16a})   implies that
$$
\zeta_{n+1,k,0}   - \zeta_{n,k-1,0}  =\zeta_{n+1,n,0}   - \zeta_{n,n-1,0}  = 1-1=0 .
$$

For $j=1$, (\ref{eq16g}) follows from (\ref{eq16d}) since it implies that
$$
\varepsilon_{n+1,k,1}   - \varepsilon_{n,k,1}  - \varepsilon_{n,k-1,1}   =
 \binom{n}{k}  -  \binom{n-1}{k}  -  \binom{n-1}{k-1} =0
$$
by (\ref{bin}).

For $j=1$, (\ref{eq16f}) follows from (\ref{eq16a}) since it implies that
$$
\zeta_{n+1,k,1}   - \zeta_{n,k,1}  - \zeta_{n,k-1,1}
   =
  \binom{n}{k+1}  -   \binom{n-1}{k+1}  -  \binom{n-1}{k} =0
$$
by (\ref{bin}).

To prove  (\ref{eq16f}) for $j\geq 2$, note that by  (\ref{eq15}), we have
\begin{eqnarray*}
&{}&
\zeta_{n+1,k,j}  - \zeta_{n,k-1,j}
\\ &{}&
=
\sum_{m=0}^{n  - k - j} j^m S(n-k-m,j) \left(   \binom{k+m}{k} -  \binom{k-1+m}{k-1}        \right)
\\ &{}&
=
\sum_{m=0}^{n  - k - j} j^m S(n-k-m,j)    \binom{k+m-1}{k}
\end{eqnarray*}
by (\ref{bin}). If $m=0$, we have $ \binom{k+m-1}{k} =0$, so that in the sum, $m$ starts at $m=1$. We set $\mu=m-1$ and obtain
$$
\zeta_{n+1,k,j}  - \zeta_{n,k-1,j}  =
j \sum_{\mu=0}^{n  - k -1 - j} j^{\mu} S(n-k-\mu-1,j)    \binom{k+\mu}{k}
= j \zeta_{n,k,j}
$$
as required. This proves (\ref{eq16f}) for $j\geq 2$.

The proof of  (\ref{eq16g})  for $j\geq 2$ is similar to our proof of (\ref{eq16f}) above                                                       and is omitted.

We finally consider (\ref{eq16i}). If $j\geq 2$, from (\ref{eq15}) we get
$$
\zeta_{n+1,0,j}   - j \zeta_{n,0,j} = \sum_{m=0}^{n - j} j^{m} S(n-m,j)
- j \sum_{m=0}^{n-1 - j} j^{m} S(n-1-m,j) .
$$
In the first sum, we take separately the term corresponding to $m=0$, which is $S(n,j)$, and in the remaining sum we set $\mu=m-1$ so that $\mu$ goes from $0$ to $n-j-1$. We see that the two remaining sums cancel out, and we obtain (\ref{eq16i}). The proof of (\ref{eq16i}) for $\varepsilon_{n+1,0,j}$ when $j\geq 2$ is similar.

If $j=1$, then (\ref{eq16i}) reads
$$
\zeta_{n+1,0,1}   =   \zeta_{n,0,1}  + S(n,1)=  \zeta_{n,0,1}  + 1, \qquad  \varepsilon_{n+1,0,1}   =   \varepsilon_{n,0,1}  + S(n,0) =  \varepsilon_{n,0,1} .
$$
Since by (\ref{eq16a}) and (\ref{eq16d})  we have
$$
 \zeta_{n,0,1} =  \binom{n-1}{1} =n-1
$$
and
$$
\varepsilon_{n,0,1} =  \binom{n-1}{0} =1,
$$
we obtain (\ref{eq16i}) for $j=1$.

If $j=0$, then (\ref{eq16i}) holds by  (\ref{eq16c})  and  (\ref{eq16e})  since
$S(n,0)=0$
and we interpret $S(n,-1)=0$.

{\bf Proof of Lemma~\ref{le1}.}
Differentiating (\ref{eq5}) and using (\ref{eq1}), we obtain, writing $E=\lambda e^{cz}$ for brevity,
\begin{equation} \label{eq17}
f^{(n+1)} = A_n'f + A_n f' + B_n'
= (A_n'+A_nE)f +( A_n (1-E)\alpha + B_n'      ) ,
\end{equation}
so that
\begin{equation} \label{eq18}
A_{n+1} = A_n' + A_n E , \qquad
B_{n+1}  =   A_n (1-E)\alpha + B_n'       .
\end{equation}
By (\ref{eq1}), we have $A_1=E=\lambda e^{cz}$ while (\ref{eq11}) and  (\ref{eq12}) state that
$$
A_1 = b_{1,1} e^{cz} = S(1,1) \lambda e^{cz} = \lambda e^{cz}
$$
since $S(1,1)=1$. This proves (\ref{eq11}) when $n=1$.

We now prove (\ref{eq11})  in general by induction on $n$. We assume that (\ref{eq11})  is valid for a certain $n\geq 1$. Then (\ref{eq18})  gives
\begin{eqnarray*}
A_{n+1} &=& A_n' + A_n E =
\sum_{k=1}^n b_{n,k} kc e^{kcz} + \sum_{k=1}^n b_{n,k} \lambda e^{(k+1)cz}
\end{eqnarray*}
and if we write this in the form
$$
\sum_{k=1}^{n+1} b_{n+1,k}  e^{kcz}
$$
then we must have, for $1\leq k\leq n+1$,
$$
b_{n+1,k} = kc b_{n,k}  + \lambda b_{n,k-1}
$$
with the convention that $b_{n,0}=0$.
We are assuming that for our $n$ the equation (\ref{eq12}) holds, and using it above gives
\begin{eqnarray*}
b_{n+1,k} &=& kc S(n,k) \lambda^k c^{n-k}   + \lambda  S(n,k-1) \lambda^{k-1} c^{n-k+1}
\\ &=&
\lambda^k c^{n-k+1} (k S(n,k) + S(n,k-1))
= \lambda^k c^{n-k+1} S(n+1,k)
\end{eqnarray*}
where for the last step we have used  (\ref{eq6}). Thus  (\ref{eq11}) and (\ref{eq12}) are valid with $n$ replaced by $n+1$ and this completes the induction proof for $A_n$.

From (\ref{eq1})  we obtain
$B_1= (1-E)\alpha$ while (\ref{eq13})--(\ref{eq16})  state that
$$
B_1 = \beta_{1,0} \alpha =
(\zeta_{1,0,0} - c\, \varepsilon_{1,0,0}  + (\zeta_{1,0,1} - c\, \varepsilon_{1,0,1}) c^{-1} E ) \alpha .
$$
These two expressions for $B_1$ are identical if, and only if,
$$
\zeta_{1,0,0} = 1, \,\, \varepsilon_{1,0,0} =0, \,\, \zeta_{1,0,1}=0,\,\, \varepsilon_{1,0,1}=1 ,
$$
which is true by (\ref{eq16h}).

To prove the formula for $B_n$ by induction on $n$, assume now that it is valid for a certain $n\geq 1$. We already have the formula for $A_n$. By (\ref{eq18}) and our assumptions,
\begin{eqnarray*}
B_{n+1}  &=&   A_n (1-E)\alpha + B_n'
\\ &=&
 (1-E)  \alpha  \sum_{p=1}^n b_{n,p} e^{pcz} +
 \sum_{k=0}^{n-1} \left( \sum_{j=0}^{n-k}  (\zeta_{n,k,j} - c\,  \varepsilon_{n,k,j} )  c^{n - k - j - 1} \lambda^j e^{jcz}   \right)  \alpha^{(k+1)}
\\ &+&
 \sum_{k=0}^{n-1}\alpha^{(k)} \left( \sum_{j=0}^{n-k}  (\zeta_{n,k,j} - c\,  \varepsilon_{n,k,j} )  c^{n - k - j - 1} jc \, \lambda^j e^{jcz}   \right)
\end{eqnarray*}
and if we write this in the form
$$
 \sum_{k=0}^{n} \beta_{n+1,k} \alpha^{(k)}
$$
then we must have
$$
\beta_{n+1,0} = (1-E)  \sum_{p=1}^n b_{n,p} e^{pcz}
+  \sum_{j=0}^{n}  (\zeta_{n,0,j} - c\,  \varepsilon_{n,0,j} )  c^{n  - j - 1} jc \, \lambda^j e^{jcz}    ,
$$
\begin{eqnarray}   \label{eq19b}
\beta_{n+1,k} &=&
 \sum_{j=0}^{n-k+1}  (\zeta_{n,k-1,j} - c\,  \varepsilon_{n,k-1,j} )  c^{n - k - j } \lambda^j e^{jcz}
\notag
\\ &+&
 \sum_{j=0}^{n-k}  (\zeta_{n,k,j} - c\,  \varepsilon_{n,k,j} )  c^{n - k - j - 1} jc \, \lambda^j e^{jcz}
\end{eqnarray}
for $1\leq k\leq n-1$, and
\begin{equation} \label{eq19a}
\beta_{n+1,n} =
\sum_{j=0}^{1}  (\zeta_{n,n-1,j} - c\,  \varepsilon_{n,n-1,j} )  c^{-j} \lambda^j e^{jcz} .
\end{equation}
We treat these three cases separately:

(a) We may write the expression for $\beta_{n+1,0} $  in the form
$$
 \sum_{p=1}^{n+1} \gamma_p \lambda^p e^{pcz}
$$
where for $1\leq p\leq n$ we have (again with $b_{n,0}=0$))
\begin{eqnarray*}
\gamma_p &=& b_{n,p} \lambda^{-p}
- b_{n,p-1}  \lambda^{1-p}
+  (\zeta_{n,0,p} - c\,  \varepsilon_{n,0,p} )  c^{n  - p - 1} pc
\\ &=&
S(n,p) c^{n-p} - S(n,p-1)  c^{n-p+1}
+  (\zeta_{n,0,p} - c\,  \varepsilon_{n,0,p} )  c^{n  - p - 1} pc
\end{eqnarray*}
while
$$
\gamma_{n+1} = - b_{n,n} \lambda^{-n} =-S(n,n) = -1.
$$
According to  (\ref{eq14}), we should have
\begin{equation} \label{eq19}
\beta_{n+1,0}  = \sum_{j=0}^{n+1}  (\zeta_{n+1,0,j} - c\,  \varepsilon_{n+1,0,j} )  c^{n  - j } \lambda^j e^{jcz} .
\end{equation}
This will be the case if, and only if,
\begin{equation} \label{eq20}
\zeta_{n+1,0,0}=  \varepsilon_{n+1,0,0}   = 0
\end{equation}
and for $1\leq j\leq n$,
\begin{equation} \label{eq21}
 \zeta_{n+1,0,j} - c\,  \varepsilon_{n+1,0,j}
 =
 S(n,j)   - S(n,j-1)  c
+  (\zeta_{n,0,j} - c\,  \varepsilon_{n,0,j} )  j .
\end{equation}
For $j=p=n+1$ we require that
$$
-1 =  (\zeta_{n+1,0,n+1} - c\,  \varepsilon_{n+1,0,n+1} )  c^{-1 } ,
$$
that is,
\begin{equation} \label{eq22}
\varepsilon_{n+1,0,n+1} = 1, \quad \zeta_{n+1,0,n+1}=0  ,
\end{equation}
which has been noted in the definitions in Lemma~\ref{le1}.
The equation  (\ref{eq21}) is valid if, and only if,
\begin{equation} \label{eq23}
 \zeta_{n+1,0,j} = S(n,j) + j \zeta_{n,0,j} , \quad
\varepsilon_{n+1,0,j} = S(n,j-1) + j \varepsilon_{n,0,j}  .
\end{equation}
But these are the formulas (\ref{eq16i})  that we have proved already.

(b) We next consider $\beta_{n+1,n}$, which according to (\ref{eq14}) should be equal to
$$
\sum_{j=0}^{1}  (\zeta_{n+1,n,j} - c\,  \varepsilon_{n+1,n,j} )  c^{-j} \lambda^j e^{jcz} .
$$
This agrees with  (\ref{eq19a}) if, and only if, for $j=0$ and for $j=1$ we have
$$
\zeta_{n+1,n,j} = \zeta_{n,n-1,j}, \quad \varepsilon_{n+1,n,j} = \varepsilon_{n,n-1,j} .
$$
This might  seem unlikely but it is valid since we have for all $n\geq 1$,
\begin{equation} \label{eq24}
\zeta_{n,n-1,0} =1, \,\,
\zeta_{n,n-1,1} =0,
\qquad \varepsilon_{n,n-1,0}=0, \,\,
\varepsilon_{n,n-1,1}=1 ,
\end{equation}
as noted in the definitions in Lemma~\ref{le1}.

(c) It remains to consider $\beta_{n+1,k}$
for $1\leq k\leq n-1$. According to (\ref{eq14}) we should have
$$
\beta_{n+1,k} =
\sum_{j=0}^{n+1-k}  (\zeta_{n+1,k,j} - c\,  \varepsilon_{n+1,k,j} )  c^{n - k - j } \lambda^j e^{jcz} .
$$
Comparing this to (\ref{eq19b}) we find that we need to prove that for each $j$ with
$0\leq j\leq n+1-k$, we have
$$
\zeta_{n+1,k,j} - c\,  \varepsilon_{n+1,k,j}
=
(\zeta_{n,k-1,j} - c\,  \varepsilon_{n,k-1,j} )  +
 (\zeta_{n,k,j} - c\,  \varepsilon_{n,k,j} ) j
$$
where the second term on the right is present only when $j\leq n-k$.
This means that we need to have, for $0\leq j\leq n-k$,
$$
\zeta_{n+1,k,j} = \zeta_{n,k-1,j} + j \zeta_{n,k,j} , \quad
\varepsilon_{n+1,k,j}  = \varepsilon_{n,k-1,j} + j \varepsilon_{n,k,j} .
$$
But these are the formulas (\ref{eq16f}) and (\ref{eq16g}) that we have proved already. When $j=n-k+1$, we need to have
$$
\zeta_{n+1,k,n-k+1} = \zeta_{n,k-1,n-k+1}, \qquad
\varepsilon_{n+1,k,n-k+1}  = \varepsilon_{n,k-1,n-k+1} .
$$
These equations are valid (they read $0=0$ and $1=1$) since by the definitions in Lemma~\ref{le1} we have $\zeta_{n,k,n-k}=0$ and $\varepsilon_{n,k,n-k} =1$.

This completes the proof of Lemma~\ref{le1}.

\section{The differential equation for $\alpha(z)$}\label{alfa}

We first note that by Theorem~\ref{Thm B} due to Lahiri, in the case (v) we are considering, we must have
\begin{equation} \label{eq25}
L(f)=a_{n}\lambda^{n}e^{ncz} f+ (1-a_{n}\lambda^{n}e^{ncz})\alpha
\end{equation}
where the initial definition of $L(f)$ is
$$
L(f) = \sum_{j=0}^n a_j f^{(j)} .
$$
However, we also assume that $a_n\not= 0$ and by Lahiri, we must have $a_0=0$ and, if $n=2$, then $a_{1}+a_{2}c=0$, and if $n\geq 3$, then
$$\sum_{k=1}^{n}(-1)^{k-1}a_{k}\left(\frac{2a_{n-1}}{n(n-1)a_{n}}\right)^{k-1}=0.
$$

Now by Lemma~\ref{le1},
$$
L(f) = \sum_{j=1}^n a_j f^{(j)} =
\sum_{j=1}^n a_j (A_j f + B_j)
= f \sum_{j=1}^n a_j A_j + \sum_{j=1}^n a_j  B_j .
$$
In view of this, (\ref{eq25}) reads
\begin{equation} \label{eq26}
C_1 f = C_2
\end{equation}
where
\begin{equation} \label{eq27}
C_1 = -a_{n}\lambda^{n}e^{ncz} + \sum_{j=1}^n a_j A_j
\end{equation}
and
\begin{equation} \label{eq28}
C_2 =  (1-a_{n}\lambda^{n}e^{ncz})\alpha   - \sum_{j=1}^n a_j  B_j  .
\end{equation}

We next write (\ref{eq1}) in the form
\begin{equation} \label{eq29}
\left(  \exp \left(-(\lambda /c) e^{cz} \right)  f \right)' =  \exp \left(-(\lambda /c) e^{cz} \right) (1-\lambda e^{cz}) \alpha(z) .
\end{equation}
Integrating, we obtain
\begin{equation} \label{eq30}
 \exp \left(-(\lambda /c) e^{cz} \right)  f(z) = e^{-\lambda /c} f(0) + \int_0^z \exp \left(-(\lambda /c) e^{c\zeta } \right) (1-\lambda e^{c\zeta }) \alpha(\zeta) \, d\zeta .
\end{equation}
Thus, whenever $\alpha$ is known, we get $f$ in this way. Recall that the function
$(1-\lambda e^{c z }) \alpha(z)$ is entire by (\ref{eq1a}) since $f$ is entire. 

By Lemma~\ref{le1}, we have
\begin{eqnarray} \label{eq34}
 C_1 &=&
  -a_{n}\lambda^{n}e^{ncz} + \sum_{j=1}^n a_j \sum_{k=1}^j S(j,k) \lambda^k c^{j-k} e^{kcz} \notag
\\
& = &
\sum_{p=1}^{n-1} \lambda^p e^{pcz} \left(  \sum_{j=p }^{n}  a_j S(j,p) c^{j-p}      \right)  .
\end{eqnarray}

If $C_1$ (and hence $C_2$) does not vanish identically and if we assume that $\alpha$ is a small function with respect to $f$, then (\ref{eq26}) implies after some work that the Nevanlinna characteristic function of $f$ satisfies
$$
T(r,f) = K_1 r + o(r)
$$
for some positive constant $K_1$. But then  (\ref{eq26}) can be written in the form
$$
\alpha = (D_0)^{-1} \left( C_1 f -   \sum_{k=1}^{n-1} D_k \alpha^{(k)}/\alpha \right)
$$
where each $D_k$ is a linear combination of functions of the form $e^{pcz}$ with constant coefficients. This implies that $m(r,\alpha)=K_2r+o(r)$ for some $K_2>0$. Since $N(r,\alpha) \leq N(r,1/(1-\lambda e^{cz}))=O(r)$, we obtain that $T(r,\alpha)=Kr+o(r)$ for some $K>0$, which is a contradiction with the assumption that $\alpha$ is a small function with respect to $f$. We conclude that $C_1$ is identically zero, so that by  (\ref{eq34}) we have
\begin{equation} \label{eq35a}
 \sum_{j=p }^{n}  a_j S(j,p) c^{j-p}    = 0
\end{equation}
for each $p$ with $1\leq p\leq n-1$. This can be used to solve for each $a_j$ for $1\leq j\leq n-1$ in terms of $a_n$.

When $p=n-1$,  (\ref{eq35a}) yields
$$
 a_{n-1} S(n-1,n-1) + a_n S(n,n-1) c=0 .
$$
Now $S(n-1,n-1)=1$ and $S(n,n-1)=\binom{n}{2}=n(n-1)/2$, so that this reads
$$
a_{n-1} + \frac{ n(n-1)}{2} ca_n=0.
$$
When $n=2$, this gives $a_1+ca_2=0$, which Lahiri obtained already.

When $p=n-2\geq 1$,  (\ref{eq35a}) yields
$$
 a_{n-2} S(n-2,n-2) + a_{n-1} S(n-1,n-2) c + a_n S(n,n-2) c^2=0 .
$$
Substituting here $a_{n-1} =- \frac{ n(n-1)}{2} ca_n$ we get
$$
-a_{n-2} = a_n c^2 ( -\frac{ (n-1)(n-2)}{2} \frac{ n(n-1)}{2}  + S(n,n-2) )
$$
so that, as a calculation shows,
\begin{equation} \label{eq35b}
a_{n-2} = a_n c^2 \frac{ n(n-1)(n-2)(3n-1)      } {    24    } .
\end{equation}

After these examples, we look for $a_{k}$ in the form $a_{k} = (-1)^{n-k} a_n c^{n-k} d_{k}$ for $k\geq 1$. We set $d_n=1$ so that this equation is valid also when $k=n$. It will turn out that $d_k>0$ for all $k$. We will show  that the $d_k$ are  given by the absolute values of the Stirling numbers of the first kind. First note that   (\ref{eq35a}) reads (since $S(p,p)=1$)
\begin{equation} \label{eq35c}
d_{p}    = (-1)^{n-p+1} S(n,p) - \sum_{j=p+1 }^{n-1}  (-1)^{j-p} S(j,p) d_{j} .
\end{equation}
When $p=n-1$, the sum is empty, and we get $d_{n-1} = S(n,n-1)$. The equation (\ref{eq35c}) together with $d_n=1$ determines the numbers $d_k$, $1\leq k\leq n$, uniquely for each fixed $n\geq 1$. If we can show that the numbers $|s(n,k)|$ satisfy the same equations in place of $d_k$, it will follow that $ d_k = |s(n,k)|$. Thus $d_n=s(n,n)=|s(n,n)|$. With $d_k$ replaced by $|s(n,k)|=(-1)^{n-k}s(n,k)$, (\ref{eq35c}) reads
$$
(-1)^{n-p}s(n,p)    = (-1)^{n-p+1} S(n,p) - \sum_{j=p+1 }^{n-1}  (-1)^{j-p} S(j,p) (-1)^{n-j}s(n,j)
$$
or
$$
 s(n,p)    = - S(n,p) - \sum_{j=p+1 }^{n-1}   S(j,p)  s(n,j) .
$$
On the other hand, (\ref{eq35d}) gives
$$
 \delta_{n,p} = \sum_{j=p}^n s(n,j) S(j,p) = \sum_{j=p+1}^{n-1} s(n,j) S(j,p)
 + s(n,p)S(p,p) + s(n,n) S(n,p) .
$$
Using $S(p,p)=s(n,n)=1$ and noting that $p\not= n$, so $\delta_{n,p} =0$, we get
$$
s(n,p) = - S(n,p) -  \sum_{j=p+1}^{n-1} s(n,j) S(j,p),
$$
as required. This proves that for each fixed $n\geq 2$, we have $d_k=|s(n,k)|=(-1)^{n-k}s(n,k)$ for $1\leq k\leq n$.

Thus
\begin{equation} \label{eq35e}
a_j =  (-1)^{n-j} a_n c^{n-j} d_{j} =  (-1)^{n-j} a_n c^{n-j} (-1)^{n-j}s(n,j) =
a_n c^{n-j} s(n,j)
\end{equation}
for $1\leq j\leq n$. This proves (\ref{eq35ee}) of Theorem~\ref{th0}.

Now that we have found the formulas for the numbers $a_j$, we consider the constraint obtained by Lahiri for $n\geq 3$, that is,
$$
\sum_{k=1}^{n}(-1)^{k-1}a_{k}\left(\frac{2a_{n-1}}{n(n-1)a_{n}}\right)^{k-1}=0 .
$$
We have
$$
\frac{2a_{n-1}}{n(n-1)a_{n}} = -c
$$
so we get
$$
\sum_{k=1}^{n}  c^{k-1}  a_{k}   =0
$$
and hence
$$
\sum_{k=1}^{n}   (-1)^{n-k} a_n c^{n-k} d_{k}  c^{k-1} =0
$$
and so
\begin{equation} \label{eq35aa}
\sum_{k=1}^{n}   (-1)^{k} d_{k}  =0 .
\end{equation}
When $n=2$, $d_2=1$ and $d_1=\frac{ 2(2-1)}{2} =1$ so that (\ref{eq35aa}) holds. When $n=3$, $d_3=1$, $d_2=\frac{ 3(3-1)}{2} =3$, and by (\ref{eq35b}), $d_1=
\frac{ 3(3-1)(3-2)(3\cdot 3-1)      } {    24    }= 2$, so that again (\ref{eq35aa}) holds.

In view of our formula for $d_k$, the equation  (\ref{eq35aa}) reads
$$
\sum_{k=1}^{n}   (-1)^{k} (-1)^{n-k}s(n,k)  =0,
$$
or
$$
\sum_{k=1}^{n}    s(n,k)  =0,
$$
which follows for $n\geq 2$ from (\ref{eq35aa})  by setting $t=1$ and noting that
$s(n,0)=0$.

We have in general $a_1=a_nc^{n-1}s(n,1)$ and $a_2=a_nc^{n-2}s(n,2)$ so that if $a_1a_2\not= 0$, then
$a_1/a_2=c s(n,1)/s(n,2) = -c| s(n,1)/s(n,2)  |$. When $n=2$, we have $s(2,2)=1$ and $|s(2,1)|=(2-1)!=1$, so that $a_1/a_2=-c$ as obtained by Lahiri, also when viewed in this way. For $n\geq 3$ we would not have $a_1/a_2=-c$.

The equation $C_2=0$, where $C_2$ is as in (\ref{eq28}), is a differential equation for $\alpha$ of order $n-1$. We may  now write (\ref{eq28}) in the form
\begin{equation} \label{eq35f}
0=  (1-a_{n}\lambda^{n}e^{ncz})\alpha   - \sum_{j=1}^n a_n c^{n-j} s(n,j)  \sum_{k=0}^{j-1} \beta_{j,k} \alpha^{(k)}
\end{equation}
and further
\begin{eqnarray} \label{eq35g}
0 &=&  (1-a_{n}\lambda^{n}e^{ncz})\alpha \notag
\\
&{}&  - \sum_{j=1}^n a_n c^{n-j} s(n,j)  \sum_{k=0}^{j-1}  \alpha^{(k)}   \sum_{p=0}^{j-k}  (\zeta_{j,k,p} - c\,  \varepsilon_{j,k,p} )  c^{j - k - p - 1} \lambda^p e^{pcz}
\notag
\\
&=&
 (1-a_{n}\lambda^{n}e^{ncz})\alpha  - c^{n-1} a_n \sum_{k=0}^{n-1} X_k  \alpha^{(k)}
\end{eqnarray}
where
$$
X_k = \sum_{j=k+1}^n  \sum_{p=0}^{j-k}  s(n,j) c^{-k-p}  (\zeta_{j,k,p} - c\,  \varepsilon_{j,k,p} )  \lambda^p e^{pcz}
.
$$
We can write
$$
X_k = c^{-k} \sum_{p=0}^{n-k} (\lambda/c)^p e^{pcz}  \sum_{j=\max\{k+1, p+k\}}^n
 s(n,j)  (\zeta_{j,k,p} - c\,  \varepsilon_{j,k,p} ) .
$$
This can be expressed as
\begin{eqnarray} \label{z000}
X_k & = &
c^{-k}  \sum_{j= k+1}^n
 s(n,j)  (\zeta_{j,k,0} - c\,  \varepsilon_{j,k,0} )
\notag  \\
& + &
c^{-k} \sum_{p=1}^{n-k} (\lambda/c)^p e^{pcz}  \sum_{j=p+k}^n
 s(n,j)  (\zeta_{j,k,p} - c\,  \varepsilon_{j,k,p} ) .
\end{eqnarray}
In view of the definitions in Lemma~\ref{le1}, according to which $\varepsilon_{j,k,0}=0$, $\zeta_{j,k,0}=0$ for $0\leq k\leq j-2$, and $\zeta_{j,j-1,0}=1$, we have
\begin{eqnarray} \label{z-1}
X_k & = &
c^{-k}  s(n,k+1)
\notag  \\
& + &
c^{-k} \sum_{p=1}^{n-k} (\lambda/c)^p e^{pcz}  \sum_{j=p+k}^n
 s(n,j)  (\zeta_{j,k,p} - c\,  \varepsilon_{j,k,p} ) .
\end{eqnarray}
The next lemma, whose proof we will give after some remarks, shows how one can simplify these expressions.

\begin{lemma} \label{le3}
If $n\geq 1$, $0\leq k\leq n-1$, 
we have
\begin{equation} \label{z1}
\sum_{j=p+k}^n
s(n,j)  \zeta_{j,k,p} = s(n - p,k+1)
\end{equation}
when $0\leq p\leq n-k-1$, and
\begin{equation} \label{z2}
\sum_{j=p+k}^n
s(n,j)  \varepsilon_{j,k,p} = s(n - p,k)
\end{equation}
when $0\leq p\leq n-k$.

If $p=n-k$, then
\begin{equation} \label{z3}
\sum_{j=p+k}^n
s(n,j)  \zeta_{j,k,p} = s(n,n)  \zeta_{k+p,k,p}
= 0 .
\end{equation}
\end{lemma}

We first make some remarks. By Lemma~\ref{le3}, we have
\begin{equation} \label{z4}
X_k = c^{-k} \left( s(n,k+1) + \sum_{p=1}^{n-k} (\lambda/c)^p e^{pcz}
(s(n - p,k+1) - c s(n - p,k) ) \right)
\end{equation}
with the convention that $s(m_1,m_2)=0$ when $m_2>m_1$.
In particular, since $s(n,n)=1$, $s(n-1,n)=0$, and $s(n-1,n-1)=1$,
we have
\begin{eqnarray} \label{z4a}
X_{n-1} & = & c^{1-n} \left(  1 - \lambda e^{cz}
   \right) 
\end{eqnarray}
so that the complete term in (\ref{eq35g}) involving $\alpha^{(n-1)}$ is
\begin{equation} \label{z4b}
-a_n ( 1 -  \lambda  e^{cz}           ) \alpha^{(n-1)}  .
\end{equation}
This agrees with the term involving $\alpha'$ in (\ref{eq3}) where $n=2$ (note that with this normalization, (\ref{eq35g}) ends up being (\ref{eq3}) multiplied through by $-1$).
Also
\begin{eqnarray} \label{z5}
X_0 & = &  \sum_{p=0}^{n} (\lambda/c)^p e^{pcz}
(s(n - p,1) - c s(n - p,0) ) \\
& = &
-c (\lambda/c)^n e^{ncz} +  \sum_{p=0}^{n-1} (\lambda/c)^p e^{pcz}
s(n - p,1)
\\
& = &
-c (\lambda/c)^n e^{ncz} +  \sum_{p=0}^{n-1} (\lambda/c)^p e^{pcz}
(-1)^{n - p - 1}  (n - p - 1)!
\end{eqnarray}
since $s(n,1)=(-1)^{n-1} (n-1)!$.
Hence in (\ref{eq35g}), the complete term involving $\alpha$ but not its derivatives is
\begin{equation} \label{z6}
\left(1-
a_n \sum_{p=0}^{n-1} c^{n - p - 1} \lambda^p e^{pcz}
(-1)^{n - p - 1} (n - p - 1)!  \right) \alpha .
\end{equation}
When $n=2$, this term is equal to
$$
(1-
a_2 ( -c + \lambda e^{cz} ) )
\alpha
$$
as already given in (\ref{eq3}) (note that with this normalization, (\ref{eq35g}) ends up being (\ref{eq3}) multiplied through by $-1$).

{\bf Proof of Lemma~\ref{le3}.}
We first consider (\ref{z1}). Suppose that $k\geq 0$ and $p\geq 0$. Then (\ref{z1}) can be considered for all $n$ with $n\geq k+p+1$. We prove (\ref{z1}) by a double induction. First, for each non-negative value $\kappa$ of $k+p$, we prove (\ref{z1}) by   induction on $n$ for $n\geq \kappa +1=k+p+1$. The non-negative values of $\kappa$ are  handled by induction on this value.

We begin with the value $\kappa=k+p=0$ so that $k=p=0$. Then, for $n=\kappa+1=1$, we have $s(n-p,k+1)=s(1,1)=1$, and the left hand side of (\ref{z1}) is (since $s(1,0)=0$ and $s(1,1)=1$, and by the line after (\ref{eq16c}))
\begin{equation} \label{z7a}
 \sum_{j=0}^{1}
 s(1,j)  \zeta_{j,0,0} = s(1,0)  \zeta_{0,0,0} + s(1,1)  \zeta_{1,0,0}
 =  \zeta_{1,0,0} =1 .
\end{equation}
Thus (\ref{z1}) holds in this case. Still with $\kappa=k+p=0$, suppose that (\ref{z1}) holds for a certain $n\geq 1$. Then the  left hand side of (\ref{z1}) with $n$ replaced by $n+1$ is equal to
$$
\sum_{j=0}^{n+1}
s(n+1,j)  \zeta_{j,0,0} = \sum_{j=1}^{n+1}
s(n+1,j)  \zeta_{j,0,0}
$$
since $s(n+1,0) =0$. When $j\geq 2$, we have $\zeta_{j,0,0}=0$ by  (\ref{eq16c}). Thus we only consider the single term $s(n+1,1)  \zeta_{1,0,0}$. Now $\zeta_{1,0,0}=1$ by the line after (\ref{eq16c}). So we finally get $s(n+1,1)$. But in our situation ($k=p=0$) we have $s(n+1-p,k+1)=s(n+1,1)$ so that holds (\ref{z1}) with $n$ replaced by $n+1$. This completes the base step for the induction on $\kappa$.

Suppose now that for a certain value of $\kappa\geq 1$, we have proved (\ref{z1}) whenever $k+p\leq \kappa -1$ and $n\geq k+p+1$. We wish to prove that for all values of $k\geq 0$ and $p\geq 0$ such that $k+p=\kappa$,  (\ref{z1}) holds for all $n\geq \kappa +1=k+p+1$. We prove this by induction on $n$. So suppose that $k\geq 0$ and $p\geq 0$ are given such that $k+p=\kappa$, and begin with the base step when $n=k+p+1$.

Suppose that $n= k+p+1$. Then (\ref{z1}) reads
\begin{equation} \label{z7}
 \sum_{j=p+k}^{k+p+1}
 s(k+p+1,j)  \zeta_{j,k,p} = s(k+1,k+1) =1 .
\end{equation}
On the other hand, when $j=k+p$, we have
$$
s(k+p+1,j)  \zeta_{j,k,p} = s(k+p+1,k+p)  \zeta_{k+p,k,p} =0
$$
since by the definition in Lemma~\ref{le1}, we have $\zeta_{n,k,n-k} =0$.
When $j=k+p+1$, we have
$$
s(k+p+1,j)  \zeta_{j,k,p} =  s(k+p+1,k+p+1)  \zeta_{k+p+1,k,p} =   \zeta_{k+p+1,k,p} .
$$
By the definition (\ref{eq15}) in Lemma~\ref{le1} (only the term $m=0$ in the sum remains), we have
$$
\zeta_{k+p+1,k,p} = S(p,p)=1
$$
when $p\geq 2$. If $p=1$, then by (\ref{eq16a})
$$
\zeta_{k+p+1,k,p} = \zeta(k+2,k,1) = \binom{k+1}{k+1}  =1 .
$$
If $p=0$, then by the line after (\ref{eq16c}),
$$
\zeta_{k+p+1,k,p} = \zeta_{k+1,k,0} =1 .
$$
This proves (\ref{z1}) when $n= k+p+1$.

Still keeping our fixed $k$ and $p$, suppose that (\ref{z1})   holds for some $n$ with $n\geq k+p+1$. Then the left hand side of (\ref{z1})  with $n$ replaced by $n+1$ is equal to
$$
\sum_{j=p+k}^{n+1}
s(n+1,j)  \zeta_{j,k,p} .
$$
By the recursion formula (\ref{z0}) $s(n+1,k) = s(n,k-1) - n s(n,k)$ we have (since $s(n,n+1) =0$ and $\zeta_{p+k,k,p}=0$)
\begin{eqnarray*}
&{}&
\sum_{j=p+k}^{n+1}
s(n+1,j)  \zeta_{j,k,p} =
\sum_{j=p+k}^{n+1}
(s(n,j-1) -n s(n,j) )  \zeta_{j,k,p}
\\ &=& s(n,p+k-1)  \zeta_{p+k,k,p} + \sum_{j=p+k}^{n} s(n,j)  \zeta_{j+1,k,p}
\\ &{}&
-n \sum_{j=p+k}^{n} s(n,j)   \zeta_{j,k,p} -n s(n,n+1) \zeta_{n+1,k,p}
\\ &=&
 \sum_{j=p+k}^{n} s(n,j)  \zeta_{j+1,k,p} -n \sum_{j=p+k}^{n} s(n,j)   \zeta_{j,k,p}
\\
&=&
\sum_{j=p+k}^{n} s(n,j)  \zeta_{j+1,k,p} -n s(n-p,k+1)
\end{eqnarray*}
where we have used our assumption (\ref{z1}). The equation (\ref{eq16f}) from Lemma~\ref{le2} implies that
$\zeta_{j+1,k,p}=p\zeta(j,k,p) + \zeta(j,k-1,p)$. Thus
\begin{eqnarray*}
&{}&
\sum_{j=p+k}^{n} s(n,j)  \zeta_{j+1,k,p} =
\sum_{j=p+k}^{n} s(n,j)  \zeta(j,k-1,p)
+p \sum_{j=p+k}^{n} s(n,j)  \zeta(j,k,p)
\\ &{}&
= s(n-p,k) + p s(n-p,k+1) .
\end{eqnarray*}
The equation $\sum_{j=p+k}^{n} s(n,j)  \zeta(j,k-1,p)= s(n-p,k)$ follows from our induction assumption involving $k+p$ since here the counterpart of $k+p$ is $(k-1)+p=\kappa -1$ (regarding the term $j=p+k-1$ that may seem to be missing here, note that $\zeta(p+k-1,k-1,p)=0$ by the line after (\ref{eq15})).
Hence
\begin{eqnarray*}
&{}&
\sum_{j=p+k}^{n+1}
s(n+1,j)  \zeta_{j,k,p} =
s(n-p,k) + p s(n-p,k+1) -n s(n-p,k+1)
\\
&{}&
= s(n-p,k) -(n-p) s(n-p,k+1)
= s(n+1-p,k+1)
\end{eqnarray*}
by the recursion formula (\ref{z0}). Thus (\ref{z1}) holds, whenever $k+p=\kappa$, with $n$ replaced by $n+1$. This completes the proof of (\ref{z1}).

We prove (\ref{z2}) using induction on $\kappa=k+p\geq 0$ and on $n \geq \kappa =k+p$ in the same way as above. Assume first that $\kappa=0$ so that $k=p=0$.
Since we are assuming that $n\geq 1$, the first case for the induction on $n$ is $n=1$. Then the right hand side of (\ref{z2}) is $s(n-p,k)=s(1,0)=0$ and the left hand side of (\ref{z2}) is
$$
\sum_{j=0}^1     s(1,j)  \varepsilon_{j,0,0}  = 0
$$
by the definition of $ \varepsilon_{j,0,0} =0$ in (\ref{eq16e}) of Lemma~\ref{le1}.
Thus  (\ref{z2}) holds in this case.

Remaining with $k=p=0$, suppose that  (\ref{z2}) has been proved for some $n\geq 1$. Then, when $n$ is replaced by $n+1$, the right hand side of (\ref{z2}) is $s(n+1-p,k)=s(n+1,0)=0$ and the left hand side of (\ref{z2}) is
$$
\sum_{j=0}^{n+1}
s(n+1,j)  \varepsilon_{j,0,0} =0
$$
since $ \varepsilon_{j,0,0} =0$ for all $j$. This completes the induction step for $n$ and proves  (\ref{z2}) for all $n\geq 1$ when $k=p=0$.

Suppose then that $\kappa\geq 1$ and that for all $k\geq 0$ and $p\geq 0$ with $k+p=\kappa-1$, the equation (\ref{z2}) has been proved for all $n$ such that $n\geq 1$ and $n\geq \kappa-1$. Now fix arbitrary $k\geq 0$ and $p\geq 0$ with $k+p=\kappa$. We need to prove by induction on $n$ for $n\geq \kappa$ that  (\ref{z2}) is valid. We begin with the base step where $n=\kappa=k+p$. Then  (\ref{z2})  reads (since there is only one term in the sum)
$$
s(p+k,p+k)  \varepsilon_{p+k,k,p} = s(n - p,k) =s(k,k)=1 .
$$
Now $s(p+k,p+k) =1$ and $ \varepsilon_{p+k,k,p}= 1$ by the line after (\ref{eq16}). Hence (\ref{z2}) holds in this case.

Keeping our fixed $\kappa$, for the induction step on $n$, suppose that (\ref{z2}) holds for a certain $n$ with $n\geq \kappa$. Then (\ref{z2}) with $n$ replaced by $n+1$ reads
$$
\sum_{j=p+k}^{n+1}
s(n+1,j)  \varepsilon_{j,k,p} = s(n+1 - p,k) .
$$
Using again the recursion formula (\ref{z0}), which gives $s(n+1,j) = s(n,j-1) - n s(n,j)$, we obtain
\begin{eqnarray*}
&{}&
\sum_{j=p+k}^{n+1}
s(n+1,j)  \varepsilon_{j,k,p}
=  \sum_{j=p+k}^{n+1}
(s(n,j-1) -n s(n,j) )  \varepsilon_{j,k,p}
\\ &=& s(n,p+k-1)  \varepsilon_{p+k,k,p} + \sum_{j=p+k}^{n} s(n,j)  \varepsilon_{j+1,k,p}
\\ &{}&
-n \sum_{j=p+k}^{n} s(n,j)   \varepsilon_{j,k,p} -n s(n,n+1) \varepsilon_{n+1,k,p}
\\ &=&
s(n,p+k-1) +
 \sum_{j=p+k}^{n} s(n,j)  \varepsilon_{j+1,k,p} -n \sum_{j=p+k}^{n} s(n,j)   \varepsilon_{j,k,p}
\\
&=&
s(n,p+k-1) + \sum_{j=p+k}^{n} s(n,j)  \varepsilon_{j+1,k,p} -n s(n-p,k+1)
\end{eqnarray*}
since $ \varepsilon_{p+k,k,p}=1$. By (\ref{eq16g}) of Lemma~\ref{le2}, we have
$$
\varepsilon_{j+1,k,p}   =  p \varepsilon_{j,k,p}  + \varepsilon_{j,k-1,p}
$$
so that
\begin{eqnarray*}
&{}&
\sum_{j=p+k}^{n} s(n,j)  \varepsilon_{j+1,k,p} =
\sum_{j=p+k}^{n} s(n,j)  \varepsilon_{j,k-1,p}
+p \sum_{j=p+k}^{n} s(n,j)  \varepsilon_{j,k,p}
\\ &{}&
= \sum_{j=p+k-1}^{n} s(n,j)  \varepsilon_{j,k-1,p}
- s(n,p+k-1)  \varepsilon_{p+k-1,k-1,p}
+ p s(n-p,k+1)
\\ &{}&
= s(n-p,k) - s(n,p+k-1) + p s(n-p,k+1)
\end{eqnarray*}
since $ \varepsilon_{p+k-1,k-1,p}=1$. We conclude that
\begin{eqnarray*}
&{}&
\sum_{j=p+k}^{n+1}
s(n+1,j)  \varepsilon_{j,k,p}
\\ &{}&
  =
s(n,p+k-1) + s(n-p,k) - s(n,p+k-1)
\\ &{}&
+ p s(n-p,k+1)  -n s(n-p,k+1)
\\ &{}&
= s(n-p,k) + p s(n-p,k+1)  -n s(n-p,k+1)
\\ &{}&
= s(n-p,k) -(n-p) s(n-p,k+1)
= s(n+1 - p,k)
\end{eqnarray*}
as required, where we have used the recursion formula (\ref{z0})  again. It follows that (\ref{z2}) holds with $n$ replaced by $n+1$, and the induction proof is complete. This proves (\ref{z2}) in all cases considered.

The equation (\ref{z3}), that is, $
s(n,n)  \zeta_{k+p,k,p} = 0$, is valid since $\zeta_{k+p,k,p}=0$ by  the line after (\ref{eq15}). This completes the proof of Lemma~\ref{le3}.

In view of Lemma~\ref{le3}, the differential equation (\ref{eq35g})  for $\alpha$ can be written as
\begin{eqnarray} \label{a1}
0 & = &
\left(1-
a_n \sum_{p=0}^{n-1} c^{n - p - 1} \lambda^p e^{pcz}
(-1)^{n - p - 1} (n - p - 1)!  \right) \alpha
\\
& {} &
-  a_n \sum_{k=1}^{n-1}   \alpha^{(k)}
c^{n-1-k} \left( s(n,k+1) + \sum_{p=1}^{n-k} (\lambda/c)^p e^{pcz}
(s(n - p,k+1) - c s(n - p,k) ) \right)
 .
 \notag
\end{eqnarray}
Taking into account (\ref{z4b}), we may also write
\begin{eqnarray} \label{a2}
0 & = &
\left(1-
a_n \sum_{p=0}^{n-1} c^{n - p - 1} \lambda^p e^{pcz}
(-1)^{n - p - 1} (n - p - 1)!  \right) \alpha \notag
\\
& {} &
-  a_n \sum_{k=1}^{n-2}   \alpha^{(k)}
c^{n-1-k} \left( s(n,k+1) + \sum_{p=1}^{n-k} (\lambda/c)^p e^{pcz}
(s(n - p,k+1) - c s(n - p,k) ) \right) \notag
\\
& {} &
-a_n ( 1 -  \lambda  e^{cz}           ) \alpha^{(n-1)} ,
\end{eqnarray}
which is a linear differential equation for $\alpha$ of order $n-1$. When $n=2$, the middle line in (\ref{a2}) represents an empty sum, which we understand to be equal to zero. This completes the proof of Theorem~\ref{th0}.

\section{Solving the differential equation for $\alpha$}\label{alfa2}

The solutions to the differential equation (\ref{a2}) for $\alpha$ have already been discussed when $n=2$.

\medskip

We  proceed to considering the case $n=3$. Recall that we now have to see under what circumstances, in general, the equation (\ref{a2}) has a non-trivial meromorphic solution $\alpha$ so that $(1-\lambda e^{cz})\alpha (z)$ is entire.

Equation (\ref{a2}) now reads
\begin{eqnarray} \label{a3}
0 & = &
\left(1-
a_3 \sum_{p=0}^{2} c^{2 - p } \lambda^p e^{pcz}
(-1)^{2 - p } (2 - p )!  \right) \alpha \notag
\\
& {} &
-  a_3   \alpha'
c \left( s(3,2) + \sum_{p=1}^{2} (\lambda/c)^p e^{pcz}
(s(3 - p,2) - c s(3 - p,1) ) \right) \notag
\\
& {} &
-a_3 ( 1 -  \lambda  e^{cz}           ) \alpha''
\notag
\\ & {} &
=  \left(1-
a_3 \left(  2c^2 - c   \lambda e^{cz} + \lambda^2 e^{2cz}
\right)   \right) \alpha \notag
\\
& {} &
-  a_3   \alpha'
c \left( s(3,2) +  (\lambda/c) e^{cz} (s(2,2) - c s(2,1) )
+ (\lambda/c)^2 e^{2cz}
(s(1,2) - c s(1,1) )
 \right) \notag
\\
& {} &
-a_3 ( 1 -  \lambda  e^{cz}           ) \alpha'' \notag
\\
& {} &
=  \left(1-
a_3 \left(  2c^2 - c   \lambda e^{cz} + \lambda^2 e^{2cz}
\right)   \right) \alpha \notag
\\
& {} &
-  a_3
 \left( -3 c +  \lambda e^{cz} (1 + c   )
-  \lambda^2 e^{2cz}
 \right)  \alpha'     \notag
\\
& {} &
-a_3 ( 1 -  \lambda  e^{cz}           ) \alpha'' .
\end{eqnarray}
Retaining only the last form for clarity, we have
\begin{eqnarray} \label{a4}
0  & = &  \left(1-
a_3 \left(  2c^2 - c   \lambda e^{cz} + \lambda^2 e^{2cz}
\right)   \right) \alpha \notag
\\
& {} &
-  a_3
 \left( -3 c +  \lambda e^{cz} (1 + c   )
-  \lambda^2 e^{2cz}
 \right)  \alpha'     \notag
\\
& {} &
-a_3 ( 1 -  \lambda  e^{cz}           ) \alpha''   .
\end{eqnarray}

\medskip

We may now write our equation for $\alpha$ in the form
\begin{equation}\label{alpha1}
\alpha ''+\frac{-3c+(1+c)\lambda e^{cz}-\lambda^{2}e^{2cz}}{1-\lambda e^{cz}}\alpha '-\frac{1/a_{3}-2c^{2}+c\lambda e^{cz}-\lambda^{2}e^{2cz}}{1-\lambda e^{cz}}\alpha =0.
\end{equation}
We write this briefly as $\alpha ''+a_{1}(z)\alpha '+a_{0}(z) \alpha =0$, where
$$a_{1}(z)=\frac{-3c+(1+c)\lambda e^{cz}-\lambda^{2}e^{2cz}}{1-\lambda e^{cz}}$$
and
$$a_{0}(z)=-\frac{1/a_{3}-2c^{2}+c\lambda e^{cz}-\lambda^{2}e^{2cz}}{1-\lambda e^{cz}}.$$

As is well-known, see, e.g., \cite{L}, p.~74, we may transform equation (\ref{alpha1}) into the form $$B''+A(z)B=0,$$ where
$$A(z)=a_{0}(z)-\frac{1}{4}a_{1}(z)^2-\frac{1}{2}a_{1}'(z)$$
and
$$B(z)=\alpha (z)\exp \left(\int \left(\frac{1}{2}a_{1}(z) \right) \right)=(\lambda e^{cz}-1)\alpha (z)e^{-3cz/2}e^{\frac{\lambda}{2c}e^{cz}}.$$

By elementary computation, we obtain
\begin{equation}\label{B1}
B''(z)+\left(-1-\frac{c^{2}}{4}-\lambda e^{cz}-\frac{\lambda^{2}}{4}e^{2cz}+
\left( \frac{1}{a_{3}}-1 \right) \frac{1}{\lambda e^{cz}-1}\right)B(z)=0.
\end{equation}
Then, of course,
\begin{equation}\label{eqalfa}
(\lambda e^{cz}-1)\alpha (z)=e^{3cz/2}e^{-\frac{\lambda}{2c}e^{cz}}B(z)
\end{equation}
and
\begin{equation}\label{eq}
f'(z)=\lambda e^{cz}f(z)-e^{3cz/2}e^{-\frac{\lambda}{2c}e^{cz}}B(z).
\end{equation}
Since $f(z)$ is entire, it follows from (\ref{eq}) that $B(z)$ has to be entire as well and so, from (\ref{eqalfa}), that $(1-\lambda e^{cz} ) \alpha(z)$ is entire.

\medskip

(1) If $a_{3}=1$, then equation (\ref{B1}) takes the form
\begin{equation}\label{B2}
B''(z)+\left(-1-\frac{c^{2}}{4}-\lambda e^{cz}-\frac{\lambda^{2}}{4}e^{2cz}\right)B(z)=0.
\end{equation}
Then it is classically well-known that all solutions to (\ref{B2}) are entire functions, forming a two-dimensional linear space. Observe that we may use Mathematica to solve (\ref{B2}):
\begin{eqnarray} \label{solB2}
B(z)
& = &
(e^{cz})^{\sqrt{1+4/c^2}/2}e^{-\frac{\lambda}{2c}e^{cz}}
\left( C_{1} U \Biggl(\frac{1}{2}(1+\sqrt{1+4/c^2})+\frac{1}{c},1+\sqrt{1+4/c^2},\frac{\lambda}{c}e^{cz} \right)+ \nonumber \\
& + &
 C_{2} L_{-(2+c+c\sqrt{1+4/c^2})/(2c)}^{\sqrt{1+4/c^2}} \left(\frac{\lambda}{c}e^{cz} \right) \Biggr), 
\end{eqnarray}
where $U(a,b,z)$, resp.\ $L_{\alpha}^{\beta}(z)$, stands for the Tricomi's confluent hypergeometric function, resp.\ the generalized Laguerre function.

\bigskip

There are explicit series expansions for the functions \hfil\break
$
U \left(\frac{1}{2}(1+\sqrt{1+4/c^2})+\frac{1}{c},1+\sqrt{1+4/c^2},\frac{\lambda}{c}e^{cz} \right)
$
and \hfil\break
$
L_{-(2+c+c\sqrt{1+4/c^2})/(2c)}^{\sqrt{1+4/c^2}} \left(\frac{\lambda}{c}e^{cz} \right) 
$, suggesting that their growth, in very rough terms, is similar to the growth of the function  $e^{\frac{\lambda}{c}e^{cz}}$. Consideration of (\ref{eq})  now shows that $f$ cannot grow essentially faster, so that $\alpha$ is not a small function of $f$. We do not attempt to give a more careful argument here.

\medskip

In  special cases, we may obtain a more explicit solution. To this end, take $2c=-3$. In this case, the square-root in (\ref{solB2}) simplifies and taking now, for simplicity, $C_{1}=1$ and $C_{2}=0$ in (\ref{solB2}) and ignoring unimportant multiplicative constants, the solution $B(z)$ reduces to
$$B(z)= ( \lambda e^{-3z/2} - 1 )e^{\frac{5z}{4}+\frac{\lambda}{3}e^{-3z/2}},$$
and $\alpha (z)$ to $\alpha (z)=  e^{-z} \exp ((2\lambda/3) e^{-3z/2})$. Solving now for $f(z)$ from (\ref{eq1}), we obtain, in terms of the incomplete gamma function $\Gamma[a,z] = \int_z^{\infty} t^{a-1} e^{-t}\, dt$,
\begin{eqnarray*}
&{}&
f(z)=e^{-(2\lambda/3) e^{-3 z/2} }  \Biggl( C + (4 \cdot 6^{1/3} (- \lambda)^{2/3} )^{-1} \times 
\\ &{}&
 \times  \left(  4 \ \Gamma[2/3, -(4/3)\lambda e^{-3 z/2} ] + 
     3 \, \Gamma[5/3,  -(4/3)\lambda e^{-3 z/2} ]  \right)  \Biggr) 
\end{eqnarray*}
where $C$ is any complex constant. 
Here it is easier to verify directly that $\alpha$ is not a small function of $f$.
We omit the details.

\medskip

(2) If $a_{3}\neq 1$, then all poles of
$$A(z)= -1-\frac{c^{2}}{4}-\lambda e^{cz}-\frac{\lambda^{2}}{4}e^{2cz}+ \left(\frac{1}{a_{3}}-1 \right)\frac{1}{\lambda e^{cz}-1}
$$
are simple (and indeed at such points $z$, where $\lambda e^{cz}=1$). Then all non-trivial meromorphic solutions to (\ref{B1}), if such solutions exist at all, are linearly dependent. Similarly, all possible solutions $\alpha (z)$ to (\ref{a4}) are linearly dependent, hence this applies to
$$(\lambda e^{cz}-1)\alpha (z)=e^{3cz/2}e^{-\frac{\lambda}{2c}e^{cz}}B(z)$$
as well. 

Writing (\ref{a4}) as a differential equation for the function $g(z) = (1 - \lambda e^{cz}  ) \alpha(z)$ that 
is required to be entire, and considering the equation at the points $z$ where 
$1 - \lambda e^{cz}=0$, we see that at such points $
(-1/a_3 +   1) g(z) = 0 
$ so that $g(z)=0$. This shows that the function $\alpha = g/ (1 - \lambda e^{cz}  ) $ is entire (if a non-trivial solution that is meromorphic in the plane exists at all for (\ref{a4})). 

Now writing (\ref{a4}) in the form
\begin{equation} \label{f2}
0= \frac{ \alpha'' } { \alpha } + \left(  ( \lambda e^{cz} - c) + \frac{ 2c } {  \lambda e^{cz} - 1}     \right) \frac{ \alpha' } { \alpha } + 
\left(  -1 -   ( \lambda e^{cz} - c) +  \frac{ A_0 } {  \lambda e^{cz} - 1}     \right) ,
\end{equation}
where $A_0 = 1/a_3 - 1 + c - 2c^2$, 
and considering it on the sub-arcs of $ \{z\colon |z|=r\} $ where $|e^{cz}|$ is large, we find by two approximate integrations that $\log \alpha'(z)$ and then $\alpha(z)$ must be of the same order of magnitude as $-   \lambda e^{cz} $ and $e^{  -(\lambda/c) e^{cz}  }$, respectively, on certain sub-arcs of these sub-arcs, the number of such sub-arcs being comparable to $r$. Hence by (\ref{eq}) $f$ cannot grow essentially faster than $\alpha$, and so $\alpha$ is not a small function of $f$. Again we do not attempt to give a more detailed proof here.

\medskip

For larger values of $n\geq 4$, the reasoning applied in the cases $n=2,3$ cannot be used directly. However, it seems reasonable to conjecture that there are entire solutions $\alpha$ and $f$ to the original problem only for some special values of the parameters involved.


\begin{thebibliography}{99}


\bibitem{Cha} Charalambides, C.A.,  \emph{Enumerative combinatorics}, Chapman \& Hall/CRC, 2002.

\bibitem{Hay} Hayman, W.K., \emph{Meromorphic functions}, Oxford University Press,  Oxford, 1964.

\bibitem{Lah} Lahiri, I., \emph{An entire function that shares a small function with its derivative and linear differential polynomial}, Comput. Method. Funct. Theory \textbf{23} (2023), 393--316.

\bibitem{L} Laine, I., \emph{Nevanlinna Theory and Complex Differential
Equations}, Walter de Gruyter, Berlin, 1993.


\bibitem{WL} Wang, J. and Laine, I., \emph{Uniqueness of entire functions and their derivatives}, Computational Methods and Function Theory \textbf{8} (2008), 327--338.

\bibitem{Z} Zhang, J., \emph{Entire functions sharing a small function with their derivatives}, Computational Methods and Function Theory \textbf{9} (2009), 379--389.



\end{thebibliography}
\end{document}